\documentstyle[11pt]{article}
\title{\huge On the Topology of Fibrations with Section\\
and Free Loop Spaces}

\author{Sadok Kallel
and Denis Sjerve\thanks{Research supported by NSERC grant
A7218}}
\date{}

\begin{document}
\maketitle


\def\la#1{\hbox to #1pc{\leftarrowfill}}
\def\ra#1{\hbox to #1pc{\rightarrowfill}}
\def\fract#1#2{\raise4pt\hbox{$ #1 \atop #2 $}}
\def\decdnar#1{\phantom{\hbox{$\scriptstyle{#1}$}}
\left\downarrow\vbox{\vskip10pt\hbox{$\scriptstyle{#1}$}}\right.}
\def\decupar#1{\phantom{\hbox{$\scriptstyle{#1}$}}
\left\uparrow\vbox{\vskip15pt\hbox{$\scriptstyle{#1}$}}\right.}

\parskip=1pc
\def\za{\vrule height6pt width4pt depth1pt}
\def\lrar{{\ra 2}}
\def\ad#1{\hbox{ad}_{#1}}

\def\map#1{\hbox{Map}_{#1}}
\def\bmap#1{\hbox{Map}^*_{#1}}
\def\map#1{\hbox{Map}_{#1}}
\def\ext{\Lambda}
\def\power{\Gamma}
\def\tensor{\otimes}
\def\dim{\hbox{dim}}

\def\bbz{{\bf Z}}
\def\bbf{{\bf F}}
\def\bbp{{\bf P}}
\def\bbr{{\bf R}}
\def\bba{{\bf A}}
\def\bbc{{\bf C}}
\def\bbr{{\bf R}}
\def\bbq{{\bf Q}}


\begin{abstract}
We relate the brace products of a fibration with section to the
differentials in its serre spectral sequence. In the particular case
of free loop fibrations, we establish a link between these
differentials and browder operations in the fiber. 
Applications and several calculations (for the particular case of
spheres and wedges of spheres) are given.
\end{abstract}


\noindent{\bf\Large\S1 Introduction}

Let $\zeta: F\fract{i}{\lrar} E\fract{\pi}{\lrar} B$ be a fibration
with a section $B\fract{s}{\lrar} E$. One of the main results of this
note asserts that the differentials on the spherical classes in the
serre spectral sequence for $\zeta$ are entirely determined by ``brace
products''in $\zeta$.

Brace products for a fibration with section were originally defined by James 
([J]).  Given  $\alpha\in\pi_p(B)$ and $\beta\in\pi_q(F)$, one can
take the Whitehead product $[ s_*(\alpha ), i_*(\beta )]$ in $\pi_{p+q-1}(E)$.
Since $\pi_*([ s_*(\alpha ), i_*(\beta )])=0$, one deduces from the long exact 
sequence in homotopy associated to $\zeta$ that $[ s_*(\alpha ), i_*(\beta )]$
must lift to a class (unique once the section is chosen)
$$\{\alpha, \beta \} \in \pi_{p+q-1}(F);$$
the so called {\sl brace product} of $\alpha$ and $\beta$.  Note that this 
class depends on the choice of section.  The brace product operation gives 
then a pairing
$$\{,\}: \pi_p(B)\times\pi_q(F)\lrar \pi_{p+q-1}(F).$$
Let $h:\pi_*(X)\lrar H_*(X;\bbz )$ denote the Hurewicz homomorphism. Our main 
result can now be stated

\noindent{\bf Theorem} 1.1:~{\sl Let $F\rightarrow E\rightarrow B$ be
a fibration with section, and assume $B$ is simply connected.  Then in
the serre spectral sequence for $E$ (with untwisted coefficients), 
the following diagram commutes
$$\matrix{
\pi_p(B)\tensor \pi_q(F)&\fract{\{,\}}{\lrar}&\pi_{p+q-1}(F)\cr
\decdnar{}&&\decdnar{h}\cr
H_p(B, H_q(F))&&H_{p+q-1}(F)\cr
\decdnar{\cong}&&\decdnar{\cong}\cr
E^2_{p,q}&\fract{d^{p}}{\ra 3}&E^{2}_{0,p+q-1}
\cr}$$}

\noindent{\sc Remarks}: Some interpretations are in order.  First the
map $\pi_p(B)\tensor \pi_q(F)\lrar H_p(B, H_q(F))$ is of course the
composite
$$\pi_p(B)\tensor\pi_q(F)\fract{h\tensor h}{\lrar} H_p(B)\tensor H_q(F)
\fract{\nu}{\lrar} H_q(F, H_p(B)),$$
where $\nu$ is a universal coefficient homomorphism. Secondly the
differential $d^p$ is really a map $E^p_{p,q}\lrar E^p_{0,p+q-1},$ but
the point is that any class in $E^2_{p,q}$ coming from
$\pi_p(B)\tensor \pi_q(F)$ actually lives until $E^p_{p,q}.$ Finally,
even though the brace product does depend on the choice of section
$s,$ commutativity of the above diagram does not (cf. \S3).

Theorem 1.1 relies in its proof on a beautiful and classical theorem
of George Whitehead [W] relating the boundary homomorphism in the
homotopy long exact sequence of a free loop fibration on a space $X$
to the Whitehead products in $X$.  More precisely, let $X$ be a finite
CW complex (based at $x_0$) and consider the evaluation fibration
$$\Omega^kX\fract{i}{\lrar}{\cal L}^kX\fract{ev}{\lrar} X\leqno{1.2}$$
where ${\cal L}^kX=\map{}(S^k, X)$ is the space of all continuous maps
from $S^k$ to $X$ (the ``$k$-th free loop space''), 
and $\Omega^kX$ is the subspace of basepoint
preserving maps. We let ${\cal L}_f^n(X)$ denote
the component containing a given map $f$.  

\noindent{\bf Theorem} 1.3: [W]~{\sl The homotopy boundary
$\partial:\pi_p(X)\rightarrow\pi_{p-1}(\Omega_f^k(X))\cong\pi_{p+k-1}(X)$
in the long exact sequence in homotopy associated to
$$\Omega_f^k(X)\lrar{\cal L}_f^k(X)\fract{ev}{\lrar} X$$
is given (up to sign) by the Whitehead product:
$\partial \alpha = [\alpha,f], \alpha\in\pi_p(X).$}

We give a sketch-proof of this theorem in \S3 and use it there to
prove 1.1.

\noindent{\sc Free Loop Spaces}: A particularly interesting application of
the theorem 1.1 occurs for the evaluation fibration 1.2 when the
connectivity of $X$ is greater than $k$. In this case 1.2 admits a
section and theorem 1.1 applies.  

In addition to loop sum, the homology ring $H_*(\Omega^k X)$ admits
a second homology operation on two variables called the browder
operation and denoted by $\lambda_k$. This operation is essential in
the calculation of the homology of iterated loop spaces (cf. [C1]
for extensive details). We quickly sketch its construction: first of
all there is an {\sl operad} map
$$\theta: S^{k-1}\times\Omega^k X \times \Omega^k X\lrar \Omega^k X$$
given as follows: a map $f\in\Omega^kX$ can be thought of as a map of
a closed unit disc in $\bbr^k$ into $X$ which sends the boundary to
basepoint. If one identifies $S^{k-1}$ with the space of pairs of
closed non-overlapping discs in $\bbr^k$, then to each pair
$(D_1,D_2)$ and to $(f,g)\in\Omega^kX$ one associates the map $\theta
(f,g)$ which is $f$ on the first disc, $g$ on the second and sends the
complement and boundary of $D_1\sqcup D_2$ to basepoint.  One then defines
$\lambda_n (x;y) := \theta_*(\iota_k, x, y)\in
H_{|x|+|y|+k-1}(\Omega^k X)$.

Let $\rho_k$ be the map
$\pi_*(X)\fract{ad_k}{\lrar}\pi_{*-k}(\Omega^kX)\fract{h}{\lrar}
H_{*-k}(\Omega^k X)$ where $ad_k$ is the adjoint isomorphism and $h$
the Hurewicz map. If we identify the spherical classes in $H_p(X)$
with classes (of the same name) in $\pi_p(X)$, then $\rho$ determines a
map from the spherical classes in $H_p(X)$ to $H_{p-k}(\Omega^kX)$.
The second main observation of this
article is

\noindent{\bf Theorem} 1.4:~{\sl Let $X$ be $k$ connected, $\beta\in
H_j(\Omega^k(X))$ and $\alpha\in H_p(X)$ two spherical classes.  Then
in the homology serre spectral sequence for $\Omega^kX\lrar
{\cal L}^kX\fract{ev}{\lrar} X$, the following relation holds
$$d^p(\alpha\tensor\beta) = \lambda_k(\rho_k(\alpha ),\beta).$$}

In the case $X=S^n$ is a sphere, $1\leq k<n$, general arguments show
that the spectral sequence collapses at $E^2$ with mod-$2$
coefficients (\S4.2).  When $n$ is odd, the same collapse occurs with
mod-$p$ coefficients.  The case $n$ even is then of greater interest
and we show the following.

Let $x\in H_{n}(S^{n})$ be the orientation class and $e\in
H_{n-k}(\Omega^kS^{n})$ be the infinite cyclic generator
representing the class of the inclusion $S^{n-k}\lrar\Omega^kS^{n}$
which is adjoint to the identity map of $S^{n}$. When $n$ even, let $a\in
H_{2n-k-1}(\Omega^kS^{n})$ be the torsion free generator (see \S4).

\noindent{\bf Corollary} 1.5:~{\sl Assume $1\leq k<n$ and $n$ is even.
Then in the homology serre spectral sequence (with integral
coefficients) for the fibration
$\Omega^kS^{n}\fract{i}{\lrar}{\cal L}^kS^{n}\fract{ev}{\lrar} S^{n}$;
$$d^{n}_{n,n-k} (x\cdot e) = 2a.$$}

\noindent{\bf Corollary} 1.6:~{\sl Suppose $1\leq k<n$ and $n$ is
even, then the Poincar\'e series for $H^*({\cal L}^kS^{n};\bbq )$ is
given as follows
$$\cases{1 + (x^n+x^{n-k})/(1-x^{2n-k-1}),& $k$ is odd\cr
(1+x^{3n-k-1})/(1-x^{n-k}),&$k$ is even.
\cr}$$}

\noindent{\bf Corollary} 1.7:~{\sl Suppose $n>2$ even and
$p$ odd.  Then in the cohomology serre spectral sequence for
${\cal L}^2S^{n}$, the mod-$p$ differentials are generated by $d_n
(x\cdot e) = x_0$, where $H^*(\Omega^2S^n;\bbz_p)$ is a tensor product
of a divided power algebra on generators $e$, $y_i$, and an exterior
algebra on generators $x_i$, $\dim (x_i) = 2(n-1)p^i-1=\dim (y_i)+1$,
$i\geq 0$.}

Corollary 1.7 has also been obtained by Fred Cohen using configuration
space model techniques (cf. [BCP]).

We can carry out similar calculations for ${\cal L}^s (W)$ where $W$ is
a bouquet of spheres. In this paper we focus on the $s=1$ case  
and there easily recover the cyclic homology description of J. Jones
and R. Cohen [C2]. More explicitly, let $W=\bigvee_{k}S^{n_i+1}$ be a
bouquet of $k$ spheres, $n_i> 1$, let $a_i\in H_{n_i+1}(W)$
be the class of the $i$-th sphere, $e_i=\rho (a_i)\in H_{n_i}(\Omega W)=
T(e_1,\ldots, e_k)$ (where $T$ is the tensor algebra).  
The map $\rho=\rho_1: H_{n_i+1}(W)\rightarrow H_{n_i}(\Omega W)$ is as defined
earlier. We prove in
\S5

\noindent{\bf Theorem} 1.8:~{\sl In the serre spectral sequence for
$\Omega W\lrar{\cal L} W\lrar W$, where $W=\bigvee_{k}S^{n_i+1}$, the
differentials are given by the cyclic operators
$$d (a_r, e_{i_1}\tensor e_{i_2}\tensor\cdots\tensor e_{i_s})
= e_r\tensor e_{i_1}\tensor\cdots\tensor
e_{i_s}-(-1)^{|e_r|(|e_{i_1}|+\cdots +|e_{i_s}|)}
e_{i_1}\tensor e_{i_2}\tensor\cdots\tensor e_{i_s}\tensor e_r$$
where $e_r=\rho (a_r)$ and $|e_j|=n_j$ is the dimension of $e_j$.}

This description yields an effective method for calculating the
homology of $W$ with mod (2) coefficients and also rational
coefficients. Following some ideas of Roos and in the case when the
spheres are equidimensional, we can show

\noindent{\bf Proposition} 1.9:~{\sl Let $W=\bigvee_{k}S^{n+1}$ and denote
by $P({\cal L} W, \bbf)$ the Poincar\'e series for $H_*({\cal L}
W,\bbf)$. Then
$$P({\cal L} W,\bbz_2 ) = 1 + (1+z)(\sum_{m\geq 1} a_mz^{mn})$$
where 
$$a_m = \sum_{d|m}{1\over d}\sum_{e|d}\mu ({d\over e})k^e =
\sum_{e|m}{1\over m}\phi ({m\over e})k^e$$
$\phi$ being the Euler $\phi$-function.}

The rational case is slightly different in the case of even spheres.

\noindent{\bf Proposition} 1.10:~{\sl For $W=\bigvee_{k}S^{n+1}$.  Then
$P({\cal L} W,\bbq ) = 1 + (1+z) (\sum_{m\geq 1} a_mz^{m n})$ where
$$
a_m = \cases{ 
\sum_{e|m}{1\over m}\phi ({m\over e})k^e, & if $n$ odd, or
$n$ even and $m$ odd\cr
\displaystyle
\sum_{{d|m}\atop d~{even}}{1\over d}\sum_{e|d}\mu ({d\over e})k^e, 
& if $n$ even and $m$ even\cr}
$$}

\noindent{\sc Remark}:
Similar results hold for the higher free loop spaces ${\cal L}^sW$ and
for this we refer to [K2].
As was pointed out to us by N. Dupont, the
above calculations recover (in particular) the following beautiful
result of Roos and Parhizgar:
$$ \dim H^{2n}({\cal L} (S^3\vee S^3);\bbq )={1\over
n}\sum_{i=1}^n2^{(i,n)}$$
where $(i,n)$ is the greatest common divisor of $i$ and $n$.

\noindent{\sc Acknowledgement}: Part of this work originated at the
CRM of Montr\'eal during the first author's visit.  The first
mentioned author would like to thank the center for its
hospitality. He would also like to thank Fred Cohen for informing him
about some unpublished work and for the proof of lemma 4.6.

\vskip 20pt


\noindent{\bf\Large\S2 Brace Products: Examples and Properties} 

\noindent{\sc Notation and Conventions}:
We often (but not always) identify a map $f : S^p\lrar X$
with its homotopy class $[f]\in\pi_p(X)$. We do so when there is no
risk of confusion and to ease notation.  We also write $ad$ for the
adjoint isomorphism
$$ad_k: \pi_{i+k}(X)\fract{\cong}{\lrar}\pi_i(\Omega^kX)$$

In the introduction we defined Brace products for a fibration
$F\fract{i}{\lrar} E\fract{p}{\lrar} B$ with a section
$B\fract{s}{\lrar}E$.  Brace products are related to Whitehead
products by the commutative diagram
$$\matrix{\pi_p(B)\tensor \pi_q(F)&\fract{\{,\}}{\lrar}&
\pi_{p+q-1}(F)\cr
\decdnar{s\tensor i}&&\decdnar{i}\cr
\pi_p(E)\tensor\pi_q(E)&\fract{[,]}{\lrar}&\pi_{p+q-1}(E)\cr}
\leqno{2.1}$$

The next three examples compute the brace product pairing for some
classes of fibrations with section.

\noindent{\sc Example} 2.2: Let $E$ be a sphere bundle over $B=S^n$
with fiber $F=S^k$ and group $O(k+1)$;
$$S^k\fract{i}{\lrar} E\fract{\pi}{\lrar} S^n.$$
This fibration is classified (up to homotopy) by a clutching function
$\mu: S^{n-1}\lrar O(k+1)$. If $E$ has
a section then the group of the bundle reduces to $O(k)$
(because the associated vector bundle does).
The map $\mu$ factors (up to homotopy)
through $S^{n-1}\lrar O(k)\hookrightarrow O(k+1),$ giving a class
$\alpha\in\pi_{n-1}O(k).$
Let $J$ be the Hopf-Whitehead construction
$$J: \pi_{n-1}(O(k))\lrar\pi_{n+k-1}(S^k).$$
Finally let $\iota_n\in\pi_n(E)$ be the class of $s: S^n\lrar E$ and 
$\iota_k\in\pi_k(E)$ be the class of the fiber. Then (up to sign)

\noindent{\bf Proposition} 2.3:~{\sl $\{\iota_n,\iota_k\}= J\alpha$ in
$\pi_{n+k-1}(S^k)$}

\noindent{\sc Proof}: We will make use of some intermediate results we
prove in \S3.  Start with the map $\alpha: S^{n-1}\lrar O(k)$. One can
think of $O(k)$ as transformations of the closed unit disc $D^k$. It
follows that $\alpha$ adjoins to a map $S^{n-1}\times D^k\lrar D^k$
and by pinching the boundary of $D^k$ we get the following commutative
diagram
$$\matrix{
S^{n-1}\times D^k&\fract{\alpha}{\lrar}&D^k\cr
\decdnar{}&&\decdnar{}\cr
S^{n-1}\wedge S^k&\fract{J}{\lrar}&S^k\cr}$$
Since $S^{n-1}\wedge S^k = S^{n+k-1}$ the bottom map can indeed be
identified with (and actually is) the $J$ homomorphism. Consider the composite
$$\phi: S^{n-1}\wedge S^k \fract{J}{\lrar} S^k\fract{i}{\lrar} E$$
We write its adjoint as a map
$g: S^{n-1}\lrar\Omega^kE$. Notice that the image of $g$ 
lies in the component
containing the fiber inclusion $i: S^k\lrar E$ and hence maps
into $\Omega^k_{\iota_k}E$.
Moreover notice that $g$ when extended to ${\cal L}^k_{\iota_k}E$
is trivial (lemma 3.4) and hence the map $g$ factors through 
the fiber $\Omega E$ of the 
inclusion $\Omega^kE\hookrightarrow{\cal L}^kE$ as follows
$$\matrix{
S^{n-1}&\fract{ad(s)}{\lrar}&\Omega E\cr
\decdnar{=}&&\downarrow\cr
S^{n-1}&\fract{g}{\lrar}&\Omega^kE&\lrar&{\cal L}^kE\cr}
\leqno{2.4}$$
The adjoint of the top map is $s: S^n\lrar E$ and the class of
this map we denote by $\iota_n$. According to lemma 3.2 we must have
that
$$\phi = [\iota_k,\iota_n]\in\pi_{n+k-1}(E)$$
Both $\phi$ and the Whitehead product map lift to $S^k$. Since the lift
of $\phi$ is $J$ and the lift of $[\iota_k,\iota_n]$ is
$\{\iota_k,\iota_n\}$, the proof is complete.

\noindent{\sc Remark}: Sphere bundles with section can be
constructed by taking a vector bundle $\zeta$ over $B=S^n$ with fiber
$F=\bbr^k$ and then compactifying fibrewise the unit disc bundle. The
new bundle (with fiber $S^{k}$) has a canonical cross section (sending
each point in $S^n$ to the point at infinity in the fiber).

\noindent{\sc Example} 2.5: It is known that the fiber of the
inclusion $X\vee X\lrar X\times X$ is $\Sigma (\Omega X\wedge \Omega
X)$ (a theorem of Ganea). Taking $X=\bbp=\bbp^{\infty}$ the infinite
complex projective space, we find that there is a fibration
$S^3\lrar\bbp\vee\bbp\lrar\bbp\times\bbp$ 
and hence after looping we obtain a fibration
$$\Omega S^3\lrar\Omega (\bbp\vee\bbp )\lrar S^1\times S^1\leqno{2.6}$$
with a section given by the composite
$$S^1\times S^1\lrar\Omega\bbp\times\Omega\bbp\fract{=}{\lrar} \Omega
(\bbp\vee *)\times\Omega (*\vee \bbp) \hookrightarrow \Omega
(\bbp\vee\bbp) \times\Omega (\bbp\vee\bbp)\fract{*}{\lrar} \Omega
(\bbp\vee\bbp).$$ 
It turns out that 2.6 has an interesting brace product given as
follows.  Denote by $a_1$ (resp. $a_2$) the generator of the second
homotopy group for the first (resp. second) copy of $\bbp$ in
$\bbp\vee\bbp$.  The fiber $S^3$ maps to $\bbp\vee\bbp$ via the
Whitehead product $[a_1,a_2]$. Taking $G:S^3\rightarrow\bbp\vee\bbp$
to be the class of the fiber, we then have that $[a_1, G] = [a_1,
[a_1,a_2]]\in\pi_4 (\bbp\vee\bbp)$ and this corresponds to the class
$\{a, G\}\in \pi_4(S^3)$ via the isomorphism $\pi_4
(\bbp\vee\bbp)\cong\pi_4(S^3)$.  It is shown in [K1]

\noindent{\bf Lemma} 2.7 [K1]: $\{a, G\}$ is the generator of
$\pi_4(S^3)=\bbz_2$.

\noindent{\sc Example} 2.8 (Saaidia): Suppose $F\lrar E\lrar B$ is a
fibration with section, and $F$ is a $G$-space with a $G$-invariant
basepoint. Consider the classifying bundle $F\lrar EG\times_GF\lrar
BG$. This fibration also admits a section and its brace products are
identified with the so-called ``secondary Eilenberg invariant'' of the
fibration $E$ (cf. [Sa]).  These invariants are fundamental in the
study of the homotopy type of the space of sections of $E$.

\noindent{\bf Brace Products and Samelson Products}

The commutator map at the level of loop spaces (better known as the
Samelson product) is related to the Whitehead product as
follows. First write $S$ for the commutator
$$S: \Omega (X)\wedge\Omega (X)\fract{}{\ra 2}\Omega (X)\ \ , \ 
(a,b)\mapsto aba^{-1}b^{-1}$$
Then the following commutes (up to sign)
$$\matrix{
\pi_{p}(\Omega X)\times\pi_{q}(\Omega X)&\fract{S}{\ra
2}&\pi_{p+q}(\Omega X)\cr 
\decdnar{ad\times ad}&&\decdnar{ad}\cr
\pi_{p+1}(X)\times\pi_{q+1}(X)&\fract{[~,~]}{\ra 2}&\pi_{p+q+1}(X).
\cr}\leqno{2.10}$$
where $ad$ is the adjoint isomorphism.  This fact (originally due to
H. Samelson) can be combined with 2.1 to show that

\noindent{\bf Lemma} 2.11:~{\sl Let $F\lrar E\lrar B$ be a fibration
with section $s$. There is a homotopy commutative diagram
$$\matrix{
\Omega B\wedge \Omega F&\fract{\{,\}}{\ra 2}&\Omega F\cr
\decdnar{\Omega s\wedge \Omega i}&&\decdnar{\Omega i}\cr
\Omega E\wedge \Omega E&\fract{S}{\ra 2}&\Omega E
\cr}$$
where the upper map (also denoted by $\{,\}$) induces James'
brace product at the level of homotopy groups}.

\noindent{\sc Proof}: The composite
$$\Omega B\wedge\Omega F\fract{S\circ (\Omega s\wedge \Omega i)}{\ra 5}
\Omega E$$
is trivial when projected into $\Omega B$ (because $\Omega p\circ\Omega i$ 
is trivial.) It then lifts to $\Omega F$ as desired. This lift is unique
up to homotopy since any two maps differ by a map
$$\Omega B\wedge\Omega F\lrar\Omega^2B\fract{\partial}{\lrar} \Omega F$$
and that this ``boundary'' map $\partial$ is null-homotopic (it is
trivial on homotopy groups because of the presence
of a section). The rest of the claim follows from 2.10.

\noindent{\bf Brace products as obstructions}

As pointed out in [J], brace products form an obstruction to
retracting the total space $E$ into the fiber $F$. They also represent
obstructions to the triviality of certain pull-back fibrations in the
postnikov tower for $B$ (see [Sa]). In what follows we exhibit yet
another obstruction expressed in terms of these brace products.

Let $F\fract{i}{\lrar} E\rightarrow B$ be a fibration of CW complexes
and consider the loop fibration
$$\Omega F\lrar\Omega E\lrar \Omega B\leqno{2.12}$$
Suppose that 2.12 has a section $s'$ and denote by * the
loop sum in $\Omega E$.  Then the composite
$$\Omega i*s': \Omega F\times\Omega B\fract{\simeq}{\ra 3}\Omega E$$
is a weak homotopy equivalence and hence an equivalence.  This
trivialization however is not necessarily an $H$-space map and its
failure to be such is measured by the commutator $(\Omega i)s'(\Omega
i)^{-1} (s')^{-1}$. We illustrate this by an example

\noindent{\sc Example} 2.13: Consider the Hopf fibering $S^1\rightarrow
S^3\rightarrow S^2$ which can be looped to a fibering
$$\Omega S^3\lrar\Omega S^2\lrar S^1$$ This has an obvious section and
as before $S^1\times\Omega S^{3}\fract{\simeq}{\lrar}
\Omega S^2$.  Notice that the left hand side is abelian (since $S^3$
is a topological group) while the right hand side $\Omega S^2$ is
not. Indeed consider the map $S^1\lrar\Omega S^2$ and
take its self commutator in $\Omega S^2$.  This commutator in homotopy
is adjoint (by the result of Samelson 2.10) to the Whitehead product 
$[\iota_2,\iota_2]=2\eta\in\pi_3(S^2)$ which is non-zero 
(here $\eta$ is the class of the hopf map). $\Omega S^2$ is hence not
abelian and the splitting $\Omega S^2\simeq S^1\times\Omega S^{3}$
is not an $H$ space splitting.

\noindent{\bf Lemma} 2.14:~{\sl 
Let $F\rightarrow E\rightarrow B$ be a fibration with section
$s$.  If the brace products in this fibration
  vanish identically, then 
$$\theta=\Omega s *\Omega i :
 \Omega B\times\Omega F\fract{\simeq}{\lrar}\Omega E$$
is an $H$-space splitting.}

\noindent{\sc Proof}: 
Here $s'=\Omega s$ and $\Omega i$ are naturally $H$-space maps
and we need only check that the following diagram homotopy commutes
$$\matrix{(\Omega B\times\Omega F)^2&~\fract{1\times\chi\times 1}{\lrar}~
(\Omega B)^2\times (\Omega F)^2~\fract{*\times *}{\ra 3}&
\Omega B\times\Omega F\cr
\decdnar{\theta^2}&&\decdnar{\theta}\cr
(\Omega E)^2&\fract{*}{\ra {13}}&\Omega E\cr}$$
where $1\times\chi\times 1$ is the shuffle map $(x,a,b,y)\mapsto (x,b,a,y)$.
Now the images of $\Omega s$ and $\Omega i$ commute in $\Omega E$
(this follows from 2.11 and from the fact that the brace products
vanish). The claim follows immediately.

\vskip 20pt


\noindent{\bf\Large\S3 Whitehead's theorem and the Proof of Theorem 1.1}

In this section we prove theorems 1.3 and 1.1 in the introduction.
Denote by $D^n$ the closed unit disc in $\bbr^n$ and by $\partial D^n
= S^{n-1}$ its boundary.  If $D^n=D^p\times D^q$, we can then write
$S^{n-1}=\partial D^n = D^p\times\partial D^q\cup\partial D^p\times
D^q$ (where the union is over $\partial D^p\times \partial D^q$). 
Let ${\cal L}^qX=\map{}(S^q,X)$ be the space of all maps from $S^q$
to $X$. We have the following pivotal lemma

\noindent{\bf Lemma} 3.1 ([W], lemma 3.3):~{\sl Start with a map
$$\phi: S^{p-1}\wedge S^{q}\lrar X$$ 
and adjoin it to get $g: S^{p-1}\lrar\Omega^q_{\alpha}X$
(where $\Omega^q_{\alpha}X$ is some component of $\Omega^qX$ containing
a representative map $\alpha$).  Suppose that $g$ extends to a map
$D^p\lrar{\cal L}^q_{\alpha}X$ and hence gives rise to an element $\beta\in
\pi_p({\cal L}^q_{\alpha}X,\Omega^q_{\alpha}X) \cong\pi_p(X)$.  Then
$$\phi = [\alpha, \beta]\in\pi_{p+q-1}(X)$$}

An alternative formulation of this lemma that is better suited to us
is as follows.

\noindent{\bf Lemma} 3.2:~{\sl Let $E$ be a space and think of $\Omega E$
as the fiber of $\Omega^qE\lrar{\cal L}^qE$. Given a composite
$$\phi: S^{p-1}\fract{\beta}{\ra 3}\Omega E\lrar\Omega^q_{\alpha}E$$
then necessarily~
$ad\phi = [\alpha, ad\beta]\in\pi_{p+q-1}E$}.

\noindent{\sc Proof}: The evaluation fibration in 1.2 extends to the
left (by looping) and we get the fibration $\Omega E\lrar
\Omega^q_{\alpha}E\lrar {\cal L}^q_{\alpha}E$. That the map $\phi:
S^{p-1}\lrar\Omega^q_{\alpha}E$ factors via $\beta$ through the fiber
$\Omega E$ is the same as having an extension diagram
$$\matrix{
S^{p-1}&\fract{\phi}{\lrar}&\Omega^q_{\alpha}E\cr
\downarrow&&\downarrow\cr
D^p&\lrar&{\cal L}^q_{\alpha}E\cr}$$
such that the element of $\pi_p({\cal L}^q_{\alpha}E,\Omega^q_{\alpha}E)
\cong\pi_p(E)$ that this diagram defines is the class of $ad\beta$.
It follows from 3.1 that $ad\phi=[\alpha, ad\beta]$.

\noindent{\bf Theorem} 3.3: [W]~{\sl The homotopy boundary
$\partial:\pi_p(X)\rightarrow\pi_{p-1}(\Omega_f^n(X))=\pi_{p+n-1}(X)$
in the long exact sequence in homotopy associated to
$$\Omega_f^q(X)\lrar{\cal L}_f^q(X)\fract{ev}{\lrar} X$$
is given (up to sign) by the Whitehead product as follows:
let $\alpha\in\pi_p(X)$, then
$\partial \alpha = ad[\alpha,f]\in\pi_{p-1}(\Omega^q_fX).$}

\noindent{\sc Proof:} Given a fibration $F\rightarrow E\rightarrow B$,
it extends to the left $\Omega B\rightarrow F$ and the boundary homomorphism
is given by the induced map in homotopy
$\pi_p(B)=\pi_{p-1}(\Omega B) \fract{\partial}{\lrar}\pi_{p-1}(F)$.
Representing $\alpha\in\pi_p(B)$ by the map of the same name,
we see that the following commutes
$$\matrix{
S^{p-1}&\fract{ad(\alpha)}{\lrar}&\Omega B\cr
\decdnar{=}&&\decdnar{}\cr
S^{p-1}&\fract{\partial\alpha}{\lrar}&F\cr}
$$
Letting $B=X$, $F=\Omega_f^qX$ and $E={\cal L}^q_f(X)$, we deduce from
lemma 3.2 that $ad^{-1}(\partial\alpha )=[\alpha, f]$ and the claim
follows.

We need one more lemma before we can proceed with the proof of
1.1. Let $\zeta: F\lrar E\lrar S^n$ be a fibration with section $s$,
and let $\mu: S^{n-1}\lrar Aut (F)$ be the cluching function. Here
$Aut(F)$ consists of based homotopy equivalences and we denote by
$\map{}^*(F,E)$ the space of based maps from $F$ into $E$. There are
inclusions
$Aut(F)\hookrightarrow\map{}^*(F,E)\hookrightarrow\map{}(F,E)$ and we
assert that

\noindent{\bf Lemma} 3.4:~{\sl There is an extension diagram
$$\matrix{
S^{n-1}&\fract{\mu}{\lrar} Aut(F)\hookrightarrow&\map{}^*(F,E)\cr
\downarrow&&\downarrow\cr
D^n&\ra{7}&\map{}(F,E)\cr}$$
such that the element $\beta\in\pi_n(\map{}(F,E),\map{}^*(F,E))\cong
\pi_n(E)$
that this defines corresponds to the class of $s:S^n\lrar E$}

\noindent{\sc Proof:}
We have the following sequence of fibrations
$F\lrar E\lrar S^n\lrar BAut(F)$ and the last map classifies the fibration
$\zeta$. By looping and letting 
$S^{n-1}\rightarrow\Omega S^n$ be the adjoint to the identity map,
we get the following diagram
$$\matrix{
&&\Omega E\cr
&&\decdnar{}\cr
S^{n-1}&\lrar&\Omega S^n&\lrar&Aut (F)\cr}$$
The lower composite, which we label $\theta$, can be identified with
the clutching map $\mu$.  If one has a section $\Omega s: \Omega
S^n\lrar \Omega E$, then $\theta$ factors through $\Omega E$ which is
the fiber of $\map{}^*(F,E)\lrar \map{}(F,E)$. The lemma follows.

\noindent{\bf Theorem} 3.5:~{\sl There is a commutative diagram
$$\matrix{
\pi_p(B)\tensor\pi_q(F)&\fract{\{,\}}{\lrar}&\pi_{p+q-1}(F)\cr
\decdnar{h\tensor h}&&\decdnar{h}\cr
H_p(B, H_q(F))&&H_{p+q-1}(F)\cr
\decdnar{\cong}&&\decdnar{\cong}\cr
E^2_{p,q}&\fract{d^p}{\lrar}&E^{2}_{0,p+q-1}
\cr}$$}

\noindent{\bf Remark} 3.6: We first explain why 3.5 is independent of
the choice of section. Suppose $F\fract{i}{\lrar} E\rightarrow B$ is
as above and assume it has two distinct sections $s_1$ and $s_2$.  Let
$\alpha\in\pi_p(B)$ and $\beta\in\pi_q(F)$.  The brace products
associated to $s_1$ and $s_2$ are given by $\{\alpha,\beta\}_1$ and
$\{\alpha,\beta\}_2$ (respectively).  Notice that
$s_1(\alpha)-s_2(\alpha )$ projects to zero in $\pi_*(B)$ and hence
must lift to a class $\alpha_F\in\pi_p(F)$. The difference element
$\{\alpha,\beta\}_1-\{\alpha,\beta\}_2$ is by definition the lift
to $\pi_*(F)$ of 
$[s_1(\alpha )-s_2(\alpha ),i_*(\beta )]=[i_*(\alpha_F),i_*(\beta )]
=i_*[\alpha_F,\beta ]\in\pi_*(E)$. It follows that
$\{\alpha,\beta\}_1-\{\alpha,\beta\}_2 = 
[\alpha_F,\beta]\in\pi_{p+q-1}(F)$. This Whitehead product in $F$
necessarily maps to zero in $H_*(F)$ by the Hurewicz homomorphism and
this is enough to show that the composite $h\circ\{,\}$ in the top
half of the diagram in 3.5 is independent of the choice of section as
asserted.

\noindent{\sc Proof of 3.5}:
Let $\alpha :S^p\to B$ represent a class in $\pi_p(B)$. Consider 
the pullback diagram
$$\matrix{F&\lrar&F\cr
\downarrow&&\downarrow\cr
E'&\lrar&E\cr
\downarrow&&\downarrow\cr
S^p&\fract{\alpha}{\lrar}&B\cr}$$
By naturality of the serre spectral sequence it suffices to prove the theorem
for the pull back fibration $F\lrar E'\lrar S^p.$ In other words we must prove
that the following diagram commutes:

$$\matrix{
\pi_p(S^p)\tensor\pi_q(F)&\fract{\{,\}}{\lrar}&\pi_{p+q-1}(F)\cr
\decdnar{h\tensor h}&&\decdnar{h}\cr
H_p(S^p, H_q(F))&&H_{p+q-1}(F)\cr
\decdnar{\cong}&&\decdnar{\cong}\cr
E^2_{p,q}&&E^{2}_{0,p+q-1}\cr
\decdnar{\cong}&&\decdnar{\cong}\cr
E^p_{p,q}&\fract{d^p}{\lrar}&E^{p}_{0,p+q-1}
\cr}$$

Now associated to $F\rightarrow E'\rightarrow S^p$ is a Wang sequence
$$\cdots\lrar H_i(F)\lrar H_i(E')\lrar H_{i-p}(F)
\fract{\tau_*}{\lrar} H_{i-1}(F)\lrar\cdots$$
where $\tau_*$ is determined in terms of the clutching
function of the bundle. Recall that this clutching function is given by
a map
$$\mu: S^{p-1}\times F\lrar F$$
whose homotopy class determines the bundle (up to fiber homotopy).
Identifying $H_{i-p}(F)$ with $E^2_{p,i-p}$ and $H_{i-1}(F)$ with
$E^2_{0,i-1}$ it isn't hard to see that $\tau_* = d^p:
E^2_{p,i-p}\rightarrow E^2_{0,i-1}$ (see [W], p:332).

Choose a basepoint $p\in F$.  Given $\beta: S^q\lrar F$ representing a
spherical class (of the same name) in $H_q(F)$, then $\tau_*$ can be
made explicit as follows.  We first have an isomorphism $H_q(F)\cong
H_{p+q}(S^p\wedge F)$ and the class $\beta$ is represented under this
isomorphism by a map $S^p\wedge S^q\rightarrow S^p\wedge F$. Writing
$D^{p+q}=D^p\times D^q$ and $\partial D^{p+q} = (D^p\times\partial
D^q)\cup (\partial D^p\times D^q)$, we can represent $\beta$
as a map of pairs
$$(D^{p+q},\partial D^{p+q})\lrar (D^p\times F, D^p\times p\cup
\partial D^p\times F).$$
The map on the second component is the boundary map $\partial$
and it can be prolonged into $F$
$$\tau: \partial D^{p+q}\fract{\partial}{\lrar} D^p\times p\cup
\partial D^p\times F\fract{}{\ra 3} F\leqno{3.7}$$
by collapsing $D^p\times p$ to $p\in F$ and sending $\partial
D^p\times F=S^{p-1}\times F$ to $F$ via the clutching function
$\mu$. (This is possible since $\mu(\partial D^p\times p)=p\in F$.)
The composite in 3.7 is a map $S^{p+q-1}\lrar F$ whose Hurewicz
image gives a class in $H_{p+q-1}(F)$. This class is exactly $\tau_*
(\beta ) = d^{p}(\beta )$.

Note at this point that the map $\tau$ gives rise by restriction to a map
$$\matrix{
\partial D^p\times D^q&\cr
\downarrow&\cr
S^{p-1}\wedge S^q&\lrar S^{p-1}\wedge F
\fract{\mu}{\lrar }F\fract{i}{\lrar} E
\cr}$$
The horizontal composite adjoins to a map $\theta: S^{p-1}\lrar
\Omega^qE$ and the component it lies in contains
the map $\beta: S^q\lrar F\lrar E$. 
By precomposing and using lemma 3.4, one gets the following 
extension diagram
$$\matrix{
S^{p-1}&\lrar&Aut(F)&\lrar&\Omega^qE\cr
\decdnar{}&&\decdnar{}&&\decdnar{}\cr
D^p&\lrar&\map{}(F,E)&\lrar&{\cal L}^qE\cr}$$
and the homotopy class this defines is given by (lemma 3.4)
$$s(S^p)\in \pi_p(E)\cong \pi_p({\cal L}^q_{\beta}E,\Omega^q_{\beta}E).$$
One can now apply lemma 3.1 directly to obtain
$$i\circ\tau = [s(\alpha ), i(\beta)]~~\hbox{in}~\pi_{p+q-1}(E).$$
Both maps lift to $F$; the LHS lifts to $\tau$ and the RHS lifts to
$\{\alpha, \beta\}:S^{p+q-1}\rightarrow F$.  Notice that in homology,
the Hurewicz images of $i_*\circ\tau_*$ and $[s(\alpha ), i(\beta)]_*$
are zero in $H_{p+q-1}(E)$ (in the first case because of the Wang
exact sequence and in the second because of a known property of
Whitehead products). It follows by the Wang exact sequence again that
the class in the image of $h\circ \{\alpha, \beta\}_*$ in
$H_{p+q-1}(F)$ is also in the image of $\tau_*$ and by the arguments
above it must follow that it is exactly $\tau_*(\beta )$. The
proposition follows.

\vskip 20pt


\noindent{\bf\Large\S4 Spaces of Free Loops}

As pointed out in the introduction, the previous results apply
particularly well to (basepoint-free) mapping spaces from spheres.
Consider again the evaluation fibration
$$\Omega^kX\fract{i}{\lrar}{\cal L}^kX\fract{ev}{\lrar} X
\leqno{4.1}$$ 
When the connectivity of $X$ is at least $k$, 4.1 admits a section
(which sends a point in $X$ to the constant loop at that point).
Below we refer by the same name to a spherical class and the homotopy
class it comes from. With $\rho_k:
\pi_*(X)\fract{ad_k}{\lrar}\pi_{*-k}(\Omega^kX)\fract{h}{\lrar}
H_{*-k}(\Omega^k X)$ as in the introduction, we prove

\noindent{\bf Theorem} 4.2:~{\sl Let $X$ be $k$ connected, $\beta\in
H_j(\Omega^k(X))$ and $\alpha\in H_p(X)$ two spherical classes.  Then
in the homology serre spectral sequence for $\Omega^kX\lrar {\cal
L}^kX\fract{ev}{\lrar} X$, the following identity holds
$$d^p(\alpha\tensor\beta) = \lambda_k(\rho_k(\alpha),\beta).$$}

\noindent{\sc Proof}: 
Suppose $M$ is $k$-connected, then the evaluation fibration admits a
section and the following commute
$$\matrix{
\pi_p(M)\tensor\pi_{j+k}(M)&\fract{[,]}{\ra {3}}&
\pi_{p+j+k-1}(M)\cr
\decdnar{1\tensor ad_k}&&\decdnar{ad_k}\cr
\pi_p(M)\tensor \pi_{j}(\Omega^k M)&\fract{\{,\}}{\lrar}&
\pi_{p+j-1}(\Omega^k M)\cr
\decdnar{h\tensor h}&&\decdnar{h}\cr
H_p(M)\tensor H_j(\Omega^k M)&\fract{d^{p}}{\ra 3}&H_{p+j-1}(\Omega^k M)
\cr}$$
The bottom half commutes because of 1.1
while the top half commutes as a result of a theorem of Hansen [H].
Notice the right vertical composite is just $\rho_k$.

Next we look at the following diagram of Fred Cohen ([C1], p:215) 
$$\matrix{
\pi_p(M)\tensor\pi_{j+k}(M)&\fract{[,]}{\ra {3}}&
\pi_{p+j+k-1}(M)\cr
\decdnar{ad_k\tensor ad_k}&&\decdnar{ad_k}\cr
\pi_{p-k}(\Omega^k M)\tensor \pi_{j}(\Omega^k M)&&
\pi_{p+j-1}(\Omega^k M)\cr
\decdnar{h\tensor h}&&\decdnar{h}\cr
H_{p-k}(\Omega^k M)\tensor H_{j}(\Omega^k M)&\fract{\lambda_k}{\lrar}&
H_{p+j-1}(\Omega^k M)\cr}$$
This diagram defines the Browder operations for spherical classes
and the proof follows by direct comparison of the above two diagrams.

\noindent{\bf Remark} 4.3: When $k=1$ and $X$ is a suspension, 
the browder operation
can be described in terms of {\bf commutators} and of the
Samelson map $\rho: \pi_*(X)\fract{ad}{\lrar}\pi_{*-1}(\Omega
X)\fract{h}{\lrar} H_{*-1}(\Omega X)$.
Let $X=S^n$ (or any suspension will do), then
according to [S] the image of a Whitehead product under
$\rho$ is a commutator in $H_*(\Omega S^n)=H_*(\Omega\Sigma S^{n-1} )
= T[e]$, where $T[e]$ is a polynomial algebra on one generator $e$ of
dimension $n-1$; i.e.
$$\rho ([x, y ])= \rho x*\rho y -(-1)^{p\cdot q} \rho x*\rho y
:= [\rho x, \rho y]$$
It then follows from theorem 4.2 that
$$d (\iota \tensor y) = [\rho (\iota), y],~~\rho (\iota )\in
H_{n-1}(\Omega S^n),~y\in H_j(\Omega S^n)$$ (here $y$ is spherical of
course). We see for instance that $d(\iota\tensor \rho (\iota ))=
[\rho (\iota), \rho (\iota)] = 
0$ if $n$ is odd, and $d(\iota\tensor \rho (\iota) ) = 2$ if $n$ is even.
This last fact generalizes to higher free loop spaces.

\noindent{\bf \S4.1 Free Loop Spaces of Spheres  ${\cal L}^kS^n, 1\leq k<n$}

When $n=1, 3$ or $7$, ${\cal L}^kS^n$ is an $H$-space (since $S^n$ is)
and so the existence of a section yields a space level splitting for
these values of $n$.  Generally and for $n$ odd, the localised sphere
$S^n_{(p)}$ at an odd prime becomes an $H$-space and hence so is
${\cal L}^kS^{n}_{(p)}$. We therefore have a space level splitting for
odd $n$ and after inverting 2.  The serre spectral sequence for 4.1
collapses for odd spheres with ${\bf Z}_p$ coefficients ($p$ odd). The
case that will preoccupy us most in this section is then when $n$ is
even.

\noindent{\bf Lemma} 4.4:~{\sl Assume $1\leq k<n$. Then
under the composite
$$\rho: \pi_{2n-1}S^n\fract{ad}{\lrar}\pi_{2n-k-1}(\Omega^kS^n)
\fract{h}{\lrar} H_{2n-k-1}(\Omega^kS^{n}),$$ 
the Whitehead square maps as follows
$$\rho([\iota_n,\iota_n])=\cases{0& $n$ is odd,\cr
2x& $n$ is even\cr}$$
(here $x$ is the infinite cyclic element in
$H_{2n-k-1}(\Omega^kS^{n};A )$, $n$ even.)}

\noindent{\sc Proof}: (sketch $n=2q$)
Write $\beta_n= [\iota_n,\iota_n]$ and let $x$ be the generator of
$H_{n-k}(\Omega^kS^n)$. Then $\rho (\beta_n) = \lambda_k (x,x)$ according
to 4.2.
When $n=2q$, $\beta_{2q}$ generates an infinite cyclic
group in $\pi_{4q-1}(S^{2q})\cong{\bf Z}\oplus\hbox{torsion}$.
It is well-known (serre) that loops on an even sphere split after
localizing at any odd prime $p$;
$$\Omega^kS^{2q}\simeq_{(p)} \Omega^{k-1}S^{2q-1}\times\Omega^kS^{4q-1}$$
Under this correspondence, it turns out that
$\beta_{2q}$ maps under $\rho$ to the 
generator in $H_{4q-k-1}(\Omega^kS^{4q-1})$ (mod (p)). Moreover
it is known that $\lambda_k (x,x)=0$ mod (2) (cf. [C1]). Putting
these together yields the result.

The following is proposition 1.6 of the introduction.

\noindent{\bf Corollary} 4.5:~{\sl Assume $1\leq k<n$ and $n$ is even.
Then in the serre spectral sequence for the fibration
$\Omega^kS^{n}\fract{i}{\lrar}{\cal L}^kS^{n}\fract{ev}{\lrar} S^{n}$,
the differential $d^{n}_{n,n-k}$ is given by multiplication by $2$ on
the torsion free generator of $H_{2n-k-1}(\Omega^kS^n)$.  In
particular, $d^{n}_{n,n-k}$ is an isomorphism with rational
coefficients.}

\noindent{\sc Proof}: The differential $d^{n}_{n,n-k}$ is determined
according to diagram 4.2 by the image of the Whitehead square under
the map $\rho$ described in 4.4. The claim now follows from lemma
4.4.  

\noindent{\bf\S4.2 Rational and Mod-$2$ Calculations}

The mod-$2$ cohomology of ${\cal L}^kS^n$, $k<n$ is completely determined
according to the following lemma

\noindent{\bf Lemma} 4.6:~{\sl The serre spectral sequence for
$\Omega^kS^n\lrar{\cal L}^kS^n\lrar S^n$ collapses with mod-$2$
coefficients whenever $k<n$}.

\noindent{\sc Proof:} (Fred Cohen) Consider the suspension $\Omega^n
E: \Omega^nS^{n+q}\rightarrow\Omega^{n+1}S^{n+q+1}$ and the following
induced map of fibrations
$$\matrix{
\Omega^nS^{n+q}&\fract{\Omega^nE}{\ra 2}&\Omega^{n+1}S^{n+q+1}\cr
\downarrow&&\downarrow\cr
{\cal L}^{n}{S^{n+q}}&\fract{{\cal L}^{n}{E}}{\ra
2}&{\cal L}^n\Omega S^{n+q+1}\cr
\downarrow&&\downarrow\cr
S^{n+q}&\fract{E}{\ra 2}&\Omega{}{S^{n+q+1}}.
\cr}$$
Since $\Omega S^{n+q+1}$ is an $H$-space, then so is
${\cal L}^n\Omega S^{n+q+1}$ and consequently we have a splitting
$${\cal L}^n\Omega S^{n+q+1}\simeq 
\Omega S^{n+q+1}\times\Omega^{n+1}S^{n+q+1}.$$
It is known (cf. [C1], pp. 228-231) that the map $\Omega^iE$ is
injective in mod-$2$ homology (for all $i$) and hence in the diagram
above both fiber and base inject in $\bbz_2$-homology.  The Lemma
follows.  

\noindent{\sc Remark}: In 5.3 below, we give an alternative derivation
of this fact in the case $k=1$.

We now use proposition 4.5 to calculate $H^*({\cal L}^kS^n)$ with rational
coefficients. We also give a complete answer mod-$p$ ($p$
odd) for the case of a two fold loop space. We make use throughout
of the following standard fact. Consider the path-loop fibration
$\Omega^kS^n\lrar P\lrar \Omega^{k-1}S^n$ for $k<n.$ Then 
$$H^*(\Omega^kS^n)=Tor^{H^*(\Omega^{k-1}S^n)}(\bbf, \bbf )\leqno{4.7}$$
This follows because the Eilenberg-Moore spectral sequence collapses at the
$E^2$ term (cf. [CM]).

\noindent{\bf Proposition} 4.8:~{\sl Let $1\leq k<n$ and suppose
$n$ even. then the Poincar\'e series for $H^*({\cal L}^kS^{n};\bbq )$ is
given as follows
$$\cases{1 + (x^n+x^{n-k})/(1-x^{2n-k-1})& , $k$ is odd\cr
(1+x^{3n-k-1})/(1-x^{n-k})& , $k$ is even.
\cr}$$}

\noindent{\sc Proof}: 
When $n$ is even, one has
$H^*(\Omega S^n)= E(e_{n-1})\tensor\bbq (a_{2n-2})$, 
where $E(e_{n-1})$ is an exterior algebra on an $n-1$
dimensional generator. It then follows (see \S4) that
$$Tor^{E (e_{n-1})}(\bbq,\bbq )=\bbq (e_{n-2}),~~\hbox{and}~~
Tor^{\bbq (a_{2n-2})}(\bbq,\bbq )=E (a_{2n-3})$$
Iterating these constructions yields
$$H^*(\Omega^kS^n;\bbq ) = \cases{\bbq (e)\tensor E (a),&$k$
even\cr 
E (e)\tensor\bbq (a),&$k$ odd\cr}$$
where $\deg (e)=n-k$ and $\deg (a)=2n-k-1.$ 
Let $\iota\in H_n(S^n )$ be the generator. Then in the serre spectral
sequence for 4.1 with $\bbq$ coefficients, the class $a$ hits $e\iota$
and this differential generates all other differentials.  When $k$ is
odd, one has (up to a unit)
$$d(a^k) = e\iota a^{k-1},~~d(ea^k) = e^2\iota a^{k-1} = 0$$
The classes that survive are $1, ea^k$ and $\tau a^k$ for $k>0$.
This establishes the first claim. When $k$ is
even, $H^*({\cal L}^kS^n;\bbq )\cong\bbq (e)[1, a\iota]$ and this leads
to the second assertion. 

\noindent{\sc Remark}: The Poincar\'e series for ${\cal L} S^n$, $n$
even; $(1 + x^n+x^{n-1}-x^{2n-2} )/(1-x^{2n-2})$, is well-known and
is given for instance in [Ro].

\noindent{\bf\S4.3 Second fold (free) loop spaces}

We now determine $H_*({\cal L}^2S^{2q+2};\bbf_p)$ with $p$ odd (the
case $p=2$ having been settled in 4.6).  So recall the description of
$\Omega^2 S^{2q+2}$ over the mod-$p$ Steenrod algebra (see [C1] or
[R] for a general discussion). We have that
$$\Omega^2 S^{2q+2}\simeq_p\Omega S^{2q+1}\times\Omega^2S^{4q+3}$$
(see the proof of 4.4), and that $H^*(\Omega^2S^{4q+3})$ is given by
$$H^*(\Omega^2S^{4q+3})={\cal L} (x_0,x_1,\ldots )\tensor
\Gamma (y_1, y_2,\ldots )$$
where $|x_i|=2(2q+1)p^i-1$ and $|y_i|=2(2q+1)p^i-2$. The action of the
Steenrod algebra is given by
$$\beta (y_i)=x_i,~~\hbox{and}~~{\cal P}^1(y_i^p)=y_{i+1}.$$

\noindent{\bf Theorem} 4.9:~{\sl In the mod (p) cohomology serre spectral
sequence for $\Omega^2 S^{2q+2}\rightarrow
{\cal L}^2S^{2q+2}\rightarrow S^{2q+2}$ we have that
$$d_{4q+1}x_0 = e\cdot\iota,$$
where $e$ is the generator of $H^{2q}(\Omega S^{2q+1})\hookrightarrow
H^{2q}(\Omega^2 S^{2q+'2}) $ in the fiber
and $\iota$ is the generator of $H^{2q+2}(S^{2q+2})$ in the base.}

\noindent{\sc Proof}:
The differential $d_{4q+1}$ is described by 4.5 and is non-trivial.
The differentials vanish on the $y$'s by dimension argument.
It follows that there are no non-zero differentials on the $x_i$'s,
$i\geq 1$ since $dx_i=d(\beta y_i)=\beta dy_i=0$. The claim
follows.

\vskip 20pt


\noindent{\bf\Large\S5 The Free Loop Space of a Bouquet of Spheres}

In this section we illustrate our techniques by calculating the
homology (with field coefficients) of ${\cal L}(\bigvee_iS^{n_i+1})$
of a finite bouquet of spheres, $n_i>0$ (similar results can be otained
for the higher free loop spaces ${\cal L}^k$; cf. [K2]).

Write $W=\bigvee_i^kS^{n_i+1}$ and consider the free loop fibration
$$\Omega W\lrar{\cal L} W\lrar W\leqno{5.1}$$ 
The image of the orientation class $[S^{n_i+1}]$ in $H_{n_i+1}(W)$
will be denoted by $a_i$ and the inclusion $S^{n_i+1}\hookrightarrow
W$ by $\iota_i$.
To the $a_i$ correspond by
adjointness the $e_i\in H_{n_i}(\Omega W)$. Observe that
$W=\Sigma (\bigvee_i^k S^{n_i})$ and so as is well-known
(Bott-Samelson)
$$H_*(W) = T(e_1,\ldots, e_k )$$
where $T(e_1,\ldots, e_k )$ is the tensor algebra on the generators
$e_i$.  An element $x\in T(e_1,\ldots, e_{k})$ is a sum of basic
monomials $e_{i_1}\tensor e_{i_2}\tensor\cdots\tensor
e_{i_r}$. Note that $x$ is not spherical in general, however iterated
commutators in the $e_i$'s are.

\noindent{\bf Lemma} 5.2 (Samelson):~{\sl 
Under the Samelson map
$\rho: \pi_*(W)\fract{ad}{\lrar}\pi_{*-1}(\Omega
W)\fract{h}{\lrar} H_{*-1}(\Omega W)$ (see 4.3),
the iterated commutator
$[e_{i_1}, [e_{i_2}, [\cdots [e_{i_{r-1}},e_{i_r}]]]$ is in the image
of the iterated Whitehead product $[\iota_{i_1}, [\iota_{i_2}, [\cdots
[\iota_{i_{r-1}},\iota_{i_r}]]]$}.

This result is also quoted in 4.3. We are now in a position to make
explicit the structure of the differentials in the serre spectral
sequence for 5.1. 

\noindent{\bf Proposition} 5.3:~{\sl In the serre spectral sequence
for 5.1, the differentials are given by
$$d_{n_i+1}(a_i\tensor x )= [e_i, x] = e_i\tensor x -
(-1)^{(n_i)|x|}x\tensor e_i$$
where again $e_i= \rho (a_i)$ in $H_{n_i}(\Omega W)$ and $x\in
H_*(\Omega W)$.}

\noindent{\sc Proof of 5.3}: The result is true for $x$ spherical
according to 4.3. Suppose now that $x=e_{i_1}\tensor
e_{i_2}\tensor\cdots\tensor e_{i_r}$ and consider the iterated
commutator $[e_{i_1}, [e_{i_2}, [\cdots [e_{i_{r-1}},e_{i_r}]]]$. This
being spherical, we get
$$d\left(a\tensor [e_{i_1}, [e_{i_2},
[\cdots [e_{i_{r-1}},e_{i_r}]]]\right) = \left[e, [e_{i_1}, [e_{i_2},
[\cdots [e_{i_{r-1}},e_{i_r}]..]\right]$$
where $e=\rho (a)$.
Writing $[e_{i_1}, [e_{i_2}, [\cdots [e_{i_{r-1}},e_{i_r}]]] =
\sum_{\tau}\pm e_{i_{\tau(1)}}e_{i_{\tau(2)}} \cdots e_{i_{\tau (r)}}$
where $\tau$ ranges over the appropriate permutations of 
$\{1,\ldots, r\}$, we can rewrite this expression as 
$$\displaystyle
\sum_{\tau}d\left(a\tensor 
e_{i_{\tau(1)}}e_{i_{\tau(2)}} \cdots e_{i_{\tau (r)}}\right)
= \sum_{\tau}\left
[e, e_{i_{\tau(1)}}e_{i_{\tau(2)}} \cdots e_{i_{\tau (r)}}\right]$$
Of course we want to show that the above summands correspond.
This is essentially forced on us by the symmetry of the situation.
We give the detailed argument for the case $r=2$ (the general case
being the same but with thiker notation). So when $r=2$
$$d(a, [e_1, e_2]) = [e, e_1e_2-e_2e_1] = ee_1e_2 -ee_2e_1 -e_1e_2e +
e_2e_1e$$
where to ease notation we choose $|e_1e_2|$ to be even and $|a|$ to be odd 
to get the appropriate signs. We stipulate $e_i\neq e_j, i\neq j$. 
We know that $d(a, [e_1, e_2])=d(a, e_1e_2) - d(a, e_2e_1)$,  and hence
one of six things must happen:
\hfill\break\hbox to .1in{\hfill}
(i) $d(a, e_1e_2) = ee_1e_2 -ee_1e_2$,~$d(a, e_2e_1) = e_1e_2e - e_2e_1e$
\hfill\break\hbox to .1in{\hfill}
(ii) $d(a, e_1e_2) = ee_1e_2 + e_2e_1e$,~$d(a, e_2e_1) = ee_2e_1 + e_1e_2e$
\hfill\break\hbox to .1in{\hfill}
(iii) $d(a, e_1e_2) = ee_1e_2-e_1e_2e = [e, e_1e_2]$, ~ $d(a, e_2e_1) = 
ee_2e_1 - e_2e_1e = [e, e_2e_1]$,
\hfill\break
the other three choices are either redundant or easily ruled out.  Of
course we need rule out (i) and (ii) to obtain (iii) for the answer.

To do this we notice generally that if $\tau$ is a permutation on $k$
letters, we can consider the bouquet $W'= \bigvee^k S^{n_{\tau
(i)}+1}$ and the (obvious) ``permutation'' map $W\lrar W'$. We get
an induced loop map $\Omega W\lrar\Omega W'$ and in turn a homology map
(which we also denote by $\tau$)
$$\tau: T(e_1,\ldots, e_k )\mapsto T(e'_1,\ldots, e'_k )=
T(e_{\tau (1)},\ldots, e_{\tau (k)})$$ 
here we have written $e_{\tau (i)}$ for 
$e'_i=\tau (e_i)$.
This map is multiplicative and induces a map of spectral sequences 
(also written $\tau$). From this we deduce
$$d(a_{\tau (i)}, e_{\tau (1)}e_{\tau (2)}\cdots e_{\tau (r)} ) 
= \tau (d(a_i, e_1e_2 \cdots e_r)) \leqno{5.4}$$ 
where by definition
$\tau (e_1\cdots e_r )=e_{\tau (1)}e_{\tau (2)}\cdots e_{\tau (r)}$.
Suppose we are in the case (i) and let $\tau$
be the transposition permuting 1 and 2 (and leaving other indexes
fixed).  We then see that $\tau d(a, e_1e_2) = \tau (ee_1e_2 -ee_1e_2)
= e'e_2e_1-e'e_2e_1$. However $d(a, e_2e_1 ) = d(a, e_{\tau (1)\tau
(2)}) = e_1e_2e' - e_2e_1e'\neq \tau d(a, e_1e_2)$. Case (i) cannot
happen.

Similary for case (ii), the same argument as above with $e=e_1$ yields
$\tau d(a_1, e_1e_2) = \tau (e_1e_1e_2 + e_2e_1e_1 ) = e_2e_2e_1 +
e_1e_2e_2\neq d(e_2, e_1e_2 )$, implying that (ii) cannot happen as
well.  Case (iii) is the only case that satisfies 5.4 as is easy to
check and the proposition follows for $r=2$. The general case $r\geq
2$ is totally analogous.  

With this description available to us, we can proceed with the
calculation of $H_*({\cal L} W)$. The following discussion is valid
with any field coefficients $\bbf$. Write $W=\Sigma X$ and let $V=
{\tilde H}_*(X;\bbf )$. The tensor algebra on $V$ corresponds to
$T(V)= T(e_1,\ldots, e_{k})$.
Consider the operator
$$\tau_m: V^{\tensor m}\lrar V^{\tensor m},\ \ \ e_{i_1}\tensor
e_{i_2}\tensor\cdots\tensor e_{i_m}\mapsto
(-1)^{n_{i_1}(n_{i_2}+\cdots + n_{i_m})} e_{i_2}\tensor\cdots\tensor
e_{i_{m}} \tensor e_{i_1}$$

The operator $\tau_m$ gives an action of the cyclic group $\bbz_m$ on
$V^{\tensor m}$ and we denote by $V^{(\tau_m)}$ the invariant subspace
under this action. Proposition 5.3 then shows that
$$H_*({\cal L} W;\bbf )\cong \bigoplus_{n\geq
0} V^{\tensor n}/_{Im(1-\tau_n)}\oplus \bigoplus_{n\geq 1}
\Sigma (V^{(\tau_n)}) \leqno{5.5}$$
where the last term is the one degree suspension of $V^{(\tau_n)}$
(compare [C2]).
Clearly $\hbox{Coker}(1-\tau_m)\subset H_*({\cal L} W)$
and the kernel of $1-\tau_m$ is a copy of ker$(1-\tau_m)$
suspended one dimension higher. Since $\dim\hbox{Coker}(1-\tau_m) =
\dim(\hbox{ker}(1-\tau_m)$), it follows that
$$P(H_*({\cal L} X)) =
1 + (1+z)P(\oplus_{m\geq 1}\hbox{ker}(1-\tau_m))\leqno{5.6} $$
where $P$ is the the mod-$\bbf$ Poincar\'e series.  In what follows we
determine $P(\oplus_{m\geq 1}Ker(1-\tau_m))$ for $\bbf=\bbq$ and
$\bbz_2$.

\noindent{\sc Definitions and Notation}: \hfill\break
$\bullet$ We denote by $\tau$ the cyclic operator $\tau
(e_{i_1}\tensor e_{i_2}\tensor \cdots\tensor\cdots e_{i_m}) =
e_{i_2}\tensor \cdots\tensor e_{i_{m}}\tensor e_{i_1}$, and by
$\tau^d$ its iterate $d$-times. It is extended to operate additively 
on all of $V^{\tensor m}$. Note that
$$\tau_m = \cases{\tau,& if $n_{i_1}$ or $n_{i_2}+\cdots + n_{i_m}$ even\cr
-\tau,& if $n_{i_1}$ and $n_{i_2}+\cdots + n_{i_m}$ odd\cr}\leqno{5.7}$$
$\bullet$ 
A word $x = e_{i_1}\tensor e_{i_2}\tensor \cdots\tensor
e_{i_m}\in V^{\tensor m}$ has period $d$ if $\tau^d (x ) =x $ and $\tau^i
(x)\neq x$ for $i<d$.  Such a word must be presented in the form of
blocks each of length $d$ and hence necessarily $d|m$. For example
$e_1e_2e_3e_1e_2e_3$ for $e_i\neq e_j$ has period $d=3$ (and $m=6$ in
this case).

\noindent{\sc The ``trick'' of Roos}:
Given a word $x$ of period $d$, consider the element
$${\bar x} = x + \tau x + \tau^2x + \cdots + \tau^{d-1}x\leqno{5.8}$$ 
Then $(1-\tau){\bar x} = x-\tau^d x = 0$.
Similarly, consider the sum
$${\bar{\bar x}} = x - \tau x + \tau^2x -\cdots(-1)^{d-1}\tau^{d-1}x$$
In this case we have
$$(1+\tau){\bar{\bar x}} =\cases{0,& if $d$ is even\cr
2x&if $d$ is odd\cr}\leqno{5.9}$$
Vice-versa, it turns out that any element in ker$(1-\tau )$ (resp.
ker$(1+\tau )$) is of the form $\bar x$ (resp. ${\bar{\bar x}}$) for
some $x$; i.e.

\noindent{\bf Lemma} 5.10:~{\sl Let $y =\sum_{\nu}e_{\nu (1)}
\tensor\cdots\tensor e_{\nu (m)}\in V^{\tensor m}\subset V$ (the sum
over some finite number of permutations $\nu$ of $\{1,\ldots,
m\}$. Then $\tau (\bar y) = \bar y$ if and only if $\bar y$ is a sum
of elements of the form ${\bar x} = x + \tau x + \tau^2x + \cdots +
\tau^{d-1}x$ for $x\in V^{\tensor m}$ and $d\geq 1$.}

\noindent{\sc Proof}: We think of $\tau$ as both an operator and a
full cyclic permutation. Clearly since $\tau (\bar x)=\bar x$, then
for any $\nu$ figuring in the expression of $\bar x$, there is a
$\nu'=\tau\circ\nu$ is also in that expression. Since
the sum is finite, there is (a smallest) $d_{\nu}\geq 1$ such that
$\nu =\tau^{d_{\nu}}\circ\nu$. The element $x=e_{\nu (1)}\tensor\cdots
e_{\nu (m)}$ has order $d_{\nu}$ and ${\bar x} = x + \tau x + \tau^2x
+ \cdots + \tau^{d_{\nu}-1}x$ is in the expression of $y$. We can then
look at $y-{\bar x}$ and proceed inductively.

Similarly if $\tau (y)=-y$, then it can be checked that $y$ is a sum
of elements of the form ${\bar{\bar x}} = x - \tau x + \tau^2x \cdots
- \tau^{d-1}x$ (here $d$ is necessarily even). (To see this one can as
a first step reduce mod-$2$ then apply the previous lemma).

Most of our forthcoming calculations are based on 5.8, 5.9 and 5.10. In
fact, let $x\in V^{\tensor m}$ be of the form $e_{i_1}\tensor
e_{i_2}\tensor\cdots\tensor e_{i_m}$ (recall $|e_i| = n_i$). 
There are two cases:

\noindent
$\bullet$ The $n_i$'s are even (i.e. the spheres are odd dimensional)
in which case $\tau_m (x) = \tau (x)$ and by
5.8, $x$ gives rise to an element $\bar x$ in the kernel of $1-\tau_m$
(any $x\in T(e_1,\ldots ,e_k)$ is necessarily periodic).
\hfill\break
$\bullet$ The $n_i$'s are not all even in which case $\tau_m
(x)=\pm\tau (x)$ and $x$ gives rise to an element in
$\hbox{ker}(1-\tau_m)$ depending on the parity of $d$ and $n_i$.

This last situation doesn't occur with mod-2 coefficients which makes
the calculations easier.

\noindent{\bf Mod-$2$ Calculations}

When $\bbf=\bbz_2$ the situation simplifies for then $\tau=\tau_m$
in all cases (see 5.7) and hence by 5.8 any
$x=e_{i_1}\tensor \cdots\tensor e_{i_m} \in V^{\tensor m}$ corresponds to an
element in the kernel of $1-\tau_m$ (namely $\bar x$).  (The same is
true when $\bbf=\bbq$ and all spheres are odd.)
Since $\bar x
= \overline{\tau x}$, $\ker (1-\tau_m )$ is in one to one
correspondence with orbits of $\tau$ acting on $V^{\tensor m}$.

\noindent{\sc Terminology}: $\tau$ acts on $T(e_1,\ldots, e_{k}) =
\bigoplus_{m\geq 1} V^{\tensor m}$ by acting on each $V^{\tensor m}$ by cyclic
permutation. An orbit consists then of a monomial $e_{i_1}\tensor
e_{i_2}\tensor\cdots\tensor e_{i_m}$ (for some $m\geq 1$) together
with all of its cyclic permutations under $\tau$. The period of the
orbit is the period of any one of its elements and the dimension of
the orbit is the homological dimension of any one of its elements.

Let $f(N)$ be the number of orbits of dimension $N$ of $\tau$ acting
on $T(e_1,\ldots, e_{k})$, and let $W=\bigvee_{k}S^{n_i+1}$ as
above. Then according to 5.6 we have

\noindent{\bf Theorem} 5.11:~{\sl 
$P(H_*({\cal L} W;\bbz_2)) = 1+ (1+z)\sum_{N\geq 1}f(N) z^{N}$}

Starting with $k$ homology classes $e_1,\ldots, e_k$, of respective
dimensions $n_1,\ldots, n_k$, and fixing an integer $N\geq 1$, we can
calculate $f(N)$ as follows. Consider all possible partitions ${\cal
P} (N)$ of $N$ by elements of $(n_{1},\ldots , n_{k})$. We write any
such partition in the form $[n_{i_1}, \ldots , n_{i_d}]$ with
$n_{i_1}+\ldots + n_{i_d}=N$. To
each partition ${\cal P} = [n_{i_1}, \ldots , n_{i_d}]$, we can let
$g({\cal P} )$ be the number of orbits made out of elements in the
corresponding tuple $(e_{i_1}, \ldots , e_{i_d})$. Then
$$f(N ) = \sum_{{\cal P}\in {\cal P} (N)}g({\cal P} )$$

\noindent{\sc Example}: 
Suppose $W=S^2\vee S^2\vee S^4\vee S^5$ and let's compute the dimension
$b_4$ of $H_4({\cal L} W;\bbz_2)$. Here
$N=4$, $n_1=1,n_2=1,n_3=3$ and $n_4=4$. We can check that we have
$8$ different partitions of $4$ by
integers taken from $\{n_1, n_2, n_3, n_4\}$; i.e
\begin{eqnarray*}
&&[n_1, n_1, n_1, n_1]\ ,\ [n_2, n_2, n_2, n_2]\ ,\ [n_1, n_1, n_1,
n_2]\ ,\ 
[n_1, n_1, n_2, n_2]\\
&&[n_1, n_2, n_2, n_2]\ ,\ [n_1, n_3]\ ,\ [n_2, n_3]\ ,\ 
[n_4]
\end{eqnarray*}
* To the partition $[n_1, n_1, n_1, n_1]$ corresponds the orbit of
$e_1\tensor e_1\tensor e_1\tensor e_1$ (of period 1),\hfill\break
* Similary to $[n_2, n_2, n_2, n_2]$ 
corresponds $e_2\tensor e_2\tensor e_2\tensor e_2$.
\hfill\break
* To $[n_1, n_2, n_2, n_2] $ corresponds only one orbit
represented by $e_1\tensor e_2\tensor e_2\tensor e_2$. 
The period here is also 4.
\hfill\break
* To $[n_1, n_1, n_1, n_2] $ 
corresponds $e_2\tensor e_1\tensor e_1\tensor e_1$.
\hfill\break
* To $[n_1, n_1, n_2, n_2]$ corresponds two orbits:
$e_1\tensor e_1\tensor e_2\tensor e_2$ and
$e_1\tensor e_2\tensor e_1\tensor e_2$.
The first has period 4 while the second has period 2.
\hfill\break
* To $[n_1, n_3]$ corresponds $e_1\tensor e_3$ (period 2).
\hfill\break
* To $[n_2, n_3]$ corresponds $e_2\tensor e_3$ (period 2).
\hfill\break
* To $[n_4]$ corresponds $e_4$ (period 1).\hfill\break
For $N=4$, there are then in total 9 orbits
and hence 9 homology classes (of degree 4). 
We need also do the same calculation for $N=3$ and there
we find 5 classes so in total
$$H_4({\cal L} W;\bbz_2)= (\bbz_2)^{14}$$

Theorem 5.11 can be made totally explicit in the case when the spheres
are all of the same dimension. The calculations there take the
following form.

\noindent{\bf Proposition} 5.12:~{\sl For $W=\bigvee_{k}S^{n}$, write
$P({\cal L} W,\bbz_2 ) = 1 + (1+z)(\sum_{m\geq 1} a_mz^{m
(n-1)})$. Then
$$a_m = \sum_{e|m}{1\over m}\phi ({m\over e})k^e$$
where $\phi$ is the Euler $\phi$-function.}

\noindent{\sc Proof}: Every element in $V^{\tensor m}$ is of degree
$N=m (n-1)$. Let
$$a_{m,d} = \hbox{number of orbits in $V^{\tensor m}$ of period $d$}$$
then $a_m=\bigoplus_{d|m}a_{m,d}$.
Let $f(d)$ be the number of monomials of period $d$. Since all
monomials in $V^{\tensor m}$ are periodic, we have that $\displaystyle
\sum_{d|m}f(d) = k^m$ and hence by the Mobius inversion formula (see
little appendix)
$$\displaystyle f(d) = \sum_{e|d}\mu ({d\over e})k^e$$
where $\mu$ is the Mobius function. It follows that
$\displaystyle
a_{d,m} = {f(d)\over d} = {1\over d}\sum_{e|d}\mu ({d\over e})k^e$.
We finally can
express $a_m$ slightly differently by using some known identities
\begin{eqnarray*}
a_m &=& \sum_{d|m}{1\over d}\sum_{e|d}\mu ({d\over e})k^e\\
&=&\sum_{e|m}k^e\sum_{h|{m\over e}}{1\over eh}\mu (h),~~h={d\over e}\\
&=&\sum_{e|m}k^e{1\over m}\left(\sum_{h|{m\over e}}{m\over eh}\mu (h)\right)
\end{eqnarray*}
The quantity in parenthesis corresponds to $\phi ({m\over e})$ according
to 5.17 below, and the proposition follows.

\noindent{\bf Remark} 5.13: When $k=1$, it is well-known that
$\sum_{e|m}{1\over m}\phi ({m\over e})$ and hence in that case $a_m=1$
for all $m$. With $\bbz_2$ coefficients, we then have
$$ \displaystyle
P(H_*({\cal L} S^n,\bbz_2))= 1+(1+z)\left(\sum_{m\geq 1} z^{m(n-1)}\right)
= 1 + (1+z)\left({1\over 1-z^{n-1}}-1\right) = {1+z^n \over 1-z^{n-1}}$$
But $P(H_*(S^n,\bbz_2))= 1+z^n$ and $P(H_*(\Omega S^n))= 
(1-z^{n-1})^{-1}$ and so we see that $H_*({\cal L} S^n,\bbz_2))\cong
H_*(S^n,\bbz_2)\tensor H_*(\Omega S^n)$ asserting that the spectral
sequence in 5.1 collapses with mod-2 coefficients when $W=S^n$ (as
asserted in 4.6). 

\noindent{\bf Mod-$\bbq$ Calculations}

Consider $W=\bigvee_{k}S^{n}$ and $n$ is odd. According to 5.7,
both actions of $\tau_m$ and of $\tau$ 
on $V^{\tensor m}$ coincide (for all $m$) and the
same argument as above shows that $P({\cal L} W,\bbq )=P({\cal L}
W,\bbz_2 )$.

We are then left with the case when $W$ is the wedge of
$k$ even dimensional spheres. Again we need determine the rank of
ker($1-\tau_m)$. Since in this case $(n-1)$ is odd
(corresponding to $n_{i_1}$ in 5.7), it follows that
$\tau_m=\tau$ when $(m-1)$ is even and $\tau_m=-\tau$ when
$(m-1)$ is odd (again by 5.7). When $(m-1)$ is even, an orbit
(of any period $d$) gives rise to an element in the kernel (cf. 5.8),
and when $(m-1)$ is odd, we get a kernel element only if $d$ is even
(cf. 5.8). That is

\noindent{\bf Proposition} 5.14:~{\sl As before 
$W=\bigvee_{k}S^{n}$ and $n$
even. Then $P({\cal L} W,\bbq ) = 1 + (1+z) (\sum_{m\geq 1}
a_mz^{m(n-1)})$ where
$$
\displaystyle a_m = \cases{ 
\sum_{{d|m}}{1\over d}\sum_{e|d}\mu ({d\over e})k^e, & if $m$ odd\cr
\sum_{{d|m}\atop d~{even}}{1\over d}\sum_{e|d}\mu ({d\over e})k^e, 
& if $m$ even\cr}
$$}

\noindent{\bf Remark} 5.15: when $k=1$, $a_m=1$ for $m$ odd and
$a_m=0$ for $m$ even (according to 5.16). In this case one regains the
calculation in 4.8
$$P(H_*({\cal L} S^n,\bbq)) = 1 + (1+z)(z^{n-1} + z^{3(n-1)} + \cdots )
= {1+ z^{n-1} + z^n - z^{2(n-1)}\over 1-z^{2(n-1)}}$$

\noindent{\sc Note Added}: The case $k>1$ with rational coefficients
is given in [K2] and the answer there is given by a suitable
desuspension of 5.14.

\noindent{\bf Appendix} (Mobius Inversion): An arithmetic function
$f:{\bf N}\lrar \bbc$ is said to be multiplicative if $f(n\cdot
m)=f(n)f(m)$ for all $n, m\in {\bf N}$. It turns out that if $f$ is
multiplicative then the function $g$ defined by $g(d) =
\sum_{e|d}f(d)$ is also multiplicative. It is possible to recover
$f(d)$ from knowledge of $g$ according to the following inversion
formula
$$g(d)= \sum_{e|d}f(e)~\Longleftrightarrow~ f(d) = \sum_{e|d}\mu
({d\over e})g(e)$$
Here $\mu (1) =1$, $\mu (n) = 0$ if $n$ has a square prime factor, and
$\mu(n) = (-1)^r$ if $n = p_1\cdots p_r$, $p_i\neq p_j$. A nice
discussion of all of this can be found in [JJ].  We simply record the
following easily established properties of the Moebius function $\mu$: 
$${1\over d}\sum_{e|d}\mu ({d\over e})=\cases{1,&$d=1$\cr
0,&otherwise\cr}\leqno{5.16}$$ 
and if $\phi (d)$ denotes the Euler $\phi$-function, then
$\sum_{e|d}\phi(e) = d$ and hence by Moebius inversion
$$\phi (d) = \sum_{e|d}{d\over e}\mu(e)\leqno{5.17}$$

\vskip 20pt
\addcontentsline{toc}{section}{Bibliography}
\bibliography{biblio}
\bibliographystyle{plain}

\vskip 10pt

\flushleft{
Sadok Kallel\hfill{Denis Sjerve}\\
Laboratoire AGAT\hfill{Dept. of Math., \#121-1884 Mathematics Road}\linebreak
Universit\'e Lille I, France\hfill{U. of British Columbia, Vancouver
V6T 1Z2} \linebreak
{\sc Email:} sadok.kallel@agat.univ-lille1.fr\hfill{sjer@math.ubc.ca}}

\end{document}